\documentclass[]{article}
\usepackage[utf8]{inputenc}
\usepackage[T1]{fontenc}

\usepackage[fleqn]{amsmath}
\usepackage{verbatim,amsthm}
\usepackage{graphicx}
\usepackage[colorinlistoftodos]{todonotes}
\usepackage[colorlinks=true, allcolors=blue]{hyperref}
\usepackage{amssymb,xcolor,amscd}
\usepackage{float}
\usepackage{enumitem}

\newcommand{\Nd}{\mathbb{N}}

\newcommand{\Un}{\mathbf{1}}

\newcommand{\eps}{\varepsilon}

\newcommand{\imply}{\Rightarrow}             %
\DeclareMathOperator{\Mean}{E}

\DeclareMathOperator*{\argmin}{argmin}
\renewcommand{\d}{\mathrm{d}}

\let\0=\varnothing
\let\1=\Un

\theoremstyle{plain}
 \newtheorem{thm}{\bfseries Theorem.}[section]

\newtheorem{cor}{\bfseries Corollary.}[section]
\newtheorem{rem}{\scshape\mdseries Remark}[section]

\theoremstyle{remark}

\numberwithin{equation}{section}

\usepackage{float}
\usepackage{listings}
\usepackage{enumitem}
\usepackage{authblk}

\title{On the Optimal Configuration of a Square Array Group Testing Algorithm}

\author[1]{Ugnė Čižikovienė}
\author[1]{Viktor Skorniakov\thanks{corresponding author; e-mail: viktor.skorniakov@mif.vu.lt}}
\affil[1]{Institute of Applied Mathematics, Faculty of Mathematics and Informatics, Vilnius University, Naugarduko 24, Vilnius LT-03225, Lithuania}

\date{}

\begin{document}

\maketitle

\begin{abstract}
  Up to date, only lower and upper bounds for the optimal configuration of a Square Array (A2) Group Testing (GT) algorithm are known. We establish exact analytical formulae and provide a couple of applications of our result. First, we compare the A2 GT scheme to several other classical GT schemes in terms of the gain per specimen attained at optimal configuration. Second, operating under objective Bayesian framework with the loss designed to attain minimum at optimal GT configuration, we suggest the preferred choice of the group size under natural minimal assumptions: the prior information regarding the prevalence suggests that grouping and application of A2 is better than individual testing. The same suggestion is provided for the Minimax strategy.
\end{abstract}

\section{Introduction}\label{s:intro}
The task of identification of infected patients in a given cohort is the frequent one. Though the plain consecutive testing of all individuals is an obvious solution, there are many other ways to approach that problem. The term Group Testing (GT) refers to the testing strategy when the testing of distinct specimens is replaced by the testing of groups of pooled specimens. It appears that the idea was first described in the famous paper of Dorfman \cite{dorfman_detection_1943}. He looked for a cost saving way to screen the U.S. soldiers for syphilis during the period of the World War II and suggested the following scheme. Instead of testing each individual blood sample, pool $N$ samples and test the group; if the group tests positive, retest each single individual; if the group tests negative, then all individuals in the group are healthy and no retesting is needed. It is clear that, when the prevalence of the disease is low, a small fraction of pools needs retesting thereby leading to significant cost savings.

Since the appearance of \cite{dorfman_detection_1943}, the idea came to stay to many other fields (quality control, informational sciences, environmental sciences, etc.) as well. In biomedical context, GT is widely applied to screen for infectious diseases like HIV, hepatitis and, most recently, COVID-19 (\cite{Xi-informativeness-1995}, \cite{wein_pooled_1996}, \cite{bilder_informative_2010}, \cite{may_pooled_2010}, \cite{stramer_nucleic_2011}, \cite{tebbs_two-stage_2013}, \cite{cov19_abdalhamid_assessment_2020}, \cite{cov19_de_wolff_evaluation_2020}, \cite{cov19_deckert_simulation_2020}, \cite{cov19_hogan_sample_2020}, \cite{cov19_lohse_pooling_2020}, \cite{cov19_pilcher_group_2020}, \cite{cov19_mutesa_pooled_2021}). It also appears to be a very useful technique in genetics (\cite{du_pooling_2006}, \cite{Cutler41}, \cite{najafi2016fundamentalDNA}, \cite{Chang-chang-Cao16DNA}).

In this paper, we focus on the Square Array (A2) GT algorithm introduced by Phatarfod and Sudbury \cite{Phatarfod-1994} and later generalized by Berger, Mandell and Subrahmanya \cite{Berger-2000}. The A2 operates as follows. Given $n^2$ specimens, one places them on $n\times n$ matrix and tests the pools defined by subsets corresponding to rows and columns. The cases lying on the intersections of rows and columns exhibiting positive responses are further retested and, assuming that the test is perfect, all infected are identified.

One of the most important characteristics of each GT algorithm is the \emph{optimal configuration}. To introduce the concept, consider an arbitrary GT scheme and let $T_N$ denote the total (random) number of tests performed over the cohort spanning $N$ individuals. The optimal configuration is the size of the cohort which minimizes function $N\mapsto \frac{\Mean T_N}{N}$, i.e., the expected number of tests per individual. Though, for a given prevalence, numerical solution of the optimal configuration is always possible, analytical formulae provide much more insights allowing, in particular, analytical comparisons of different GT schemes. For the case of the Dorfman scheme described above, the optimal configuration was established by Samuels \cite{samuels_exact_1978} quite long ago. However, for a couple of its modifications, namely, the modified Dorfman scheme \cite{sobel+groll:1959} and Sterrett scheme \cite{sterrett_detection_1957}, the analytical formulae, though conjectured, were unknown \cite{malinovsky_revisiting_2019} and established only recently \cite{cizikoviene2021couple}. Pretty much the same situation is with A2. To our best knowledge, only Hudgens and Kim \cite{hudgens_optimal_2011} addressed this problem and did not succeed in providing the final solution. Namely, in \cite{hudgens_optimal_2011} the authors obtained lower and upper bounds for the optimal A2 configuration leaving the exact formulae undiscovered. In this paper, we fill in the gap by providing exact solution similar to that obtained by Samuels \cite{samuels_exact_1978} and Skorniakov and Čižikovienė \cite{cizikoviene2021couple} for the case of classical Dorfman scheme and its modifications.

The paper is organized as follows. In Section \ref{s:res_and_apps}, we state our results and provide a couple of applications. Section \ref{s:discussion} is devoted to discussion. Finally, there are two appendices containing proofs, figures, and link to the web repository with tabular data.

\section{Results and Applications}\label{s:res_and_apps}

\subsection{Statement of the Results}\label{ss:results}
Before proceeding to the statement of results, we first formulate Binomial Testing Assumptions (BTA) which are assumed to hold in the remaining part of the paper by default.

\begin{itemize}
    \item[(BTA1)] The tested cohort consists of independent individuals. Each individual is infected with the same constant probability $p\in(0,1)$ (termed prevalence in the sequel).
    \item[(BTA2)] The test under consideration is perfect and there is no dilution effect. That is, pooling does not affect the performance of the test.
\end{itemize}

The way A2 operates was already described in the introductory Section \ref{s:intro}. However, we need an explicit expression for the average number of tests $\Mean T_N$ applied  under A2 to the cohort spanning $N=n^2,n\geq 2$, individuals. The latter was derived by Phatarfod and Sudbury \cite{Phatarfod-1994} and is equal to
\begin{equation}\label{e:mean_T_N_A2}
    2n + n^2\left(1-2q^n + q^{2n-1}\right),
\end{equation}
where $q:=1-p$. Therefore, in this parametrization, an average number of tests per person
\begin{equation}\label{e:t_q_n_A2}
    t(q,n):=\frac{1}{N}\Mean T_N \Big |_{N=n^2}=\frac{2}{n} + 1-2q^n+q^{2n-1},
\end{equation}
and the corresponding optimal configuration 
\begin{equation*}
    n_{opt}=n_{opt}(q):=\mathrm{argmin}_{n\in\{2,3,\ldots\}} t(q,n).
\end{equation*}
The existence, uniqueness and bounds on $n_{opt}$ for various values of $q\in(0,1)$ were established by Hudgens and Kim \cite{hudgens_optimal_2011}. Our reconsideration (given in Theorem \ref{t:main} and Corollary \ref{c:n_opt} below) aimed to sharpen their results (see discussion in Section \ref{s:discussion}) and provide the complete theoretical characterization of A2.

\begin{thm}\label{t:main}
    Let $g(q,n)=\frac{2}{n}-2q^n+q^{2n-1}=t(q,n)-1$. 
    \begin{itemize}
        \item[(i)] For $(q,n)$ ranging in $(1/2,1)\times(2,\infty)$, system of equations  
        \begin{equation}\label{e:cutSys}
            \begin{cases}
                    &\ln q = nq^n\left(1-\frac{q^{n-1}}{2}\right)\ln q\\
                    &n\ln q = -\frac{\left(1-\frac{q^{n-1}}{2}\right)}{\left(1-{q^{n-1}}\right)}
            \end{cases}
        \end{equation}
        has a unique solution $(q_*,n_*)\approx (0.748416, 4.453524)$. 
        \item[(ii)] For any fixed $q\in (q_*,1)$ and with respect to $n$, equation $g(q,n)=0$ admits two solutions $n_L,n_U:2<n_L<n_*<n_U<\infty$. On $(n_L,n_U)$, $n\mapsto g(q,n)$ attains values in $(-\infty,0)$ whereas on $(2,\infty)\setminus[n_L,n_U]$ it attains values in $(0,\infty)$.
        \item[(iii)] For any fixed $q\in (q_*,1)$, the region $(n_L,n_U)$ is the one where A2 is efficient, i.e., $t(q,n)<1$ for $n\in (n_L,n_U)$. In that region, there exists a unique (and, therefore, global) minimizer $n_{min}$ of $(2,\infty)\ni n\mapsto t(q,n)$. For $q\in[0.755,1)$, it is given by
        \begin{equation}\label{e:A2_n_optimal}
            n_{min} =\frac{1}{p^{\frac{2}{3}}} + \frac{1}{2p^{\frac{1}{3}}} +  0.2 + 3p^2 + t_*
        \end{equation}
        for some $t_*\in [0,1]$.
    \end{itemize}
    $(2,\infty)\ni n\mapsto t(q,n)$ also has a unique (and, therefore, global) maximizer located in the region $(n_U,\infty)$. For any fixed $q\in(0,q_*)$, A2 is never optimal, i.e., $(2,\infty)\ni n \mapsto t(q,n)$ attains values in $(1,\infty)$.
\end{thm}

\begin{cor}\label{c:n_opt}
Let $g(q,n)$ be as in Theorem \eqref{t:main}. Then $g(q,5)=0$ has a unique solution $q_{5}\approx 0.750209961$. For all $q\in(q_5,1)$, $n_{opt}(q)$ belongs to the set
\begin{equation}\label{e:set_of_optimal_vals}
    \left\{\left\lfloor \frac{1}{p^{\frac{2}{3}}} + \frac{1}{2p^{\frac{1}{3}}} + 3p^2 + 0.2\right\rfloor+i: i=0,1,2\right\}.
\end{equation}
\end{cor}

\begin{rem}\label{r:about_region_truncation}
In the Corollary \ref{c:n_opt}, the region $(q_5,1)\subset(q_*,1)$. An explanation for this truncation stems from the fact that $(q_5,1)$ is the region where practical application of A2 makes sense. To be more precise, applying Theorem \ref{t:main} for a fixed $q\in(q_*,q_5)$, we have that $t(q,n_{min})<1$ with $n_{min}=n_{min}(q)$ given by \eqref{e:A2_n_optimal}. However, this $n_{min}(q)\in(4,5)$ and $\min(t(q,4),t(q,5))>1$ when $q\in(q_5,q_*)$. We touch this question briefly in the discussion Section \ref{s:discussion} when talking about relation of our results to those of Hudgens and Kim \cite{hudgens_optimal_2011}. Though we do not provide a separate proof of this fact, technical details can be filled in after inspection of the proofs presented in the Appendix \ref{a:proofs}.
\end{rem}

\subsection{Examples of Applications}\label{ss:applications}

\subsubsection{Comparison to other GT schemes}\label{ss:example1}

In order to shed the light on to performance of A2, we compare it to several other GT schemes in terms of magnitude of optimal configuration $N_{opt}$ and gain across the range $p\in(0,0.249790)$ where application of A2 seems reasonable. In this comparative analysis, for a fixed $q=1-p$, we define the gain as 
\begin{equation*}
    G(q)=(1-t(q,N_*)).
\end{equation*}
Here, $t(q,N)=\frac{\Mean T_N}{N}$ is an average number of tests per person when the tested group size is equal to $N$ whereas $N_*$ stands for the unrounded optimal configuration, i.e., the minimizer of the continuous argument function $[1,\infty)\ni N\mapsto t(q,N)$. Defined this way, the gain multiplied by 100 has a meaning of an average number of tests saved per 100 persons in comparison to usual one by one testing. Preceded by several remarks, below comes a description of the a fore mentioned GT schemes. To distinguish between the schemes, we assign one letter abbreviations to each and, by making use of these letters, superscript all related quantities. E.g.,  $t^{(A2)}(q,N_*^{(A2)})$ denotes the value of $t(q,N_*)$ when A2 is the scheme under consideration.

\begin{rem}\label{r:reparametrization}
Recall that, in section \ref{ss:results}, we have introduced reparametrization of $N$ by equality $N=n^2$, where $n$ is the number of rows (columns) in the squared array used in the definition of A2. In this example and in all what follows after, this reparametrization remains in force: writing $t^{(A2)}(q,N)$, we factually mean $t$ given in equation \eqref{e:A2_n_optimal} with $n$ ranging continuously unless stated otherwise. When we want to emphasize reference to \eqref{e:A2_n_optimal}, $n$ is used instead of $N$. For all other schemes considered, reparametrizations of similar kind do not apply.
\end{rem}

\begin{rem}
The behaviour of quantities compared is more naturally interpreted in terms of the prevalence $p=1-q$. Therefore, $p$ but not $q$ is used in the accompanying graphs. Also, because of the same reason and for the sake of convenience, we quite often denote an argument of the function considered by $q$ and write an explicit formula in terms of $p$.
\end{rem}

\begin{rem}
In this example, we have chosen to operate on the continuous scale because it is much easier to perceive visually in comparison to the discrete one. However, keeping in a view the practical aspect, comprehensive numerical results (see Appendix \ref{app:tbl_and_figures}) are given on the discrete scale. Comparing both one can find out that the discrepancies are small.
\end{rem}

\smallskip\noindent\emph{Dorfman scheme {D}}. $\rhd$ The Dorfman scheme was described in the introductory Section \ref{s:intro}. For this scheme,
\begin{equation}\label{e:tD_N}
    t^{(D)}(q,N)= \frac{1}{N}+1-q^N
\end{equation}
and $N_{*}^{(D)}$ solves equation $1/N^2=-q^N\ln q$. Samuels \cite{samuels_exact_1978} has shown that rounded optimal configuration $N_{opt}^{(D)}$ is either $1+\lfloor p^{-1/2}\rfloor$ or $2+\lfloor p^{-1/2}\rfloor$. Hence, $N_{*}^{(D)}=\sqrt{p^{-1}}(1+o(1)),\ p\to 0+$. Approximation $N_*^{(D)}\approx\sqrt{p^{-1}}$ is accurate enough (see \cite{Skorniakov_Leipus_2021} for tabulated numerical results), and one can further show that 
\begin{equation*}
    t^{(D)}(q,N_*^{(D)})=2\sqrt{p}(1+o(1)),\ p\to 0+.
\end{equation*}
Consequently, $G^{(D)}(q)=(1-2\sqrt{p})(1+o(1)),\ p\to 0+$. $\lhd$

\smallskip\noindent\emph{Sterrett scheme {S}}. $\rhd$ Sterrett \cite{sterrett_detection_1957} suggested the following modification of the Dorfman\footnote{in fact, his suggestion was built on the already modified Dorfman scheme} scheme: one should retest initial positive pool sequentially one by one until appearance of the first positive case and then again apply pool testing to the remaining tail. If the remaining untested set tests positive, one should proceed recursively as previously until the remaining set tests negative or the whole set of individuals gets tested. Malinovsky and Albert \cite{malinovsky_revisiting_2019} conjectured that, for $p\in(0,(3-\sqrt{5})/2)$, rounded optimal configuration $N_{opt}^{(S)}$ lies in the set 
\begin{equation*}
    \left\{\lfloor\sqrt{2p^{-1}}\rfloor,\lfloor\sqrt{2p^{-1}}\rfloor+1,\lfloor\sqrt{2p^{-1}}\rfloor+2\right\}.
\end{equation*}
Skorniakov and Čižikovienė \cite{cizikoviene2021couple} affirmed their conjecture showing along the way that $N_{*}^{(S)}\in\left[\sqrt{2p^{-1}}-1,\sqrt{2p^{-1}}+1\right]$ for $p\in(0,(3-\sqrt{5})/2)$. From the latter result it then follows that
\begin{align*}
    &N_{*}^{(S)}=\sqrt{2p^{-1}}(1+o(1)),\quad t^{(S)}(q,N_*^{(S)})=\sqrt{2p}(1+o(1)),\ p\to0+.\
\end{align*}
Hence, $G^{(S)}(q)=(1-\sqrt{2p})(1+o(1)),\ p\to 0+$. $\lhd$

\begin{figure}[ht!]
\centering
\includegraphics[width=15cm, height=12cm]{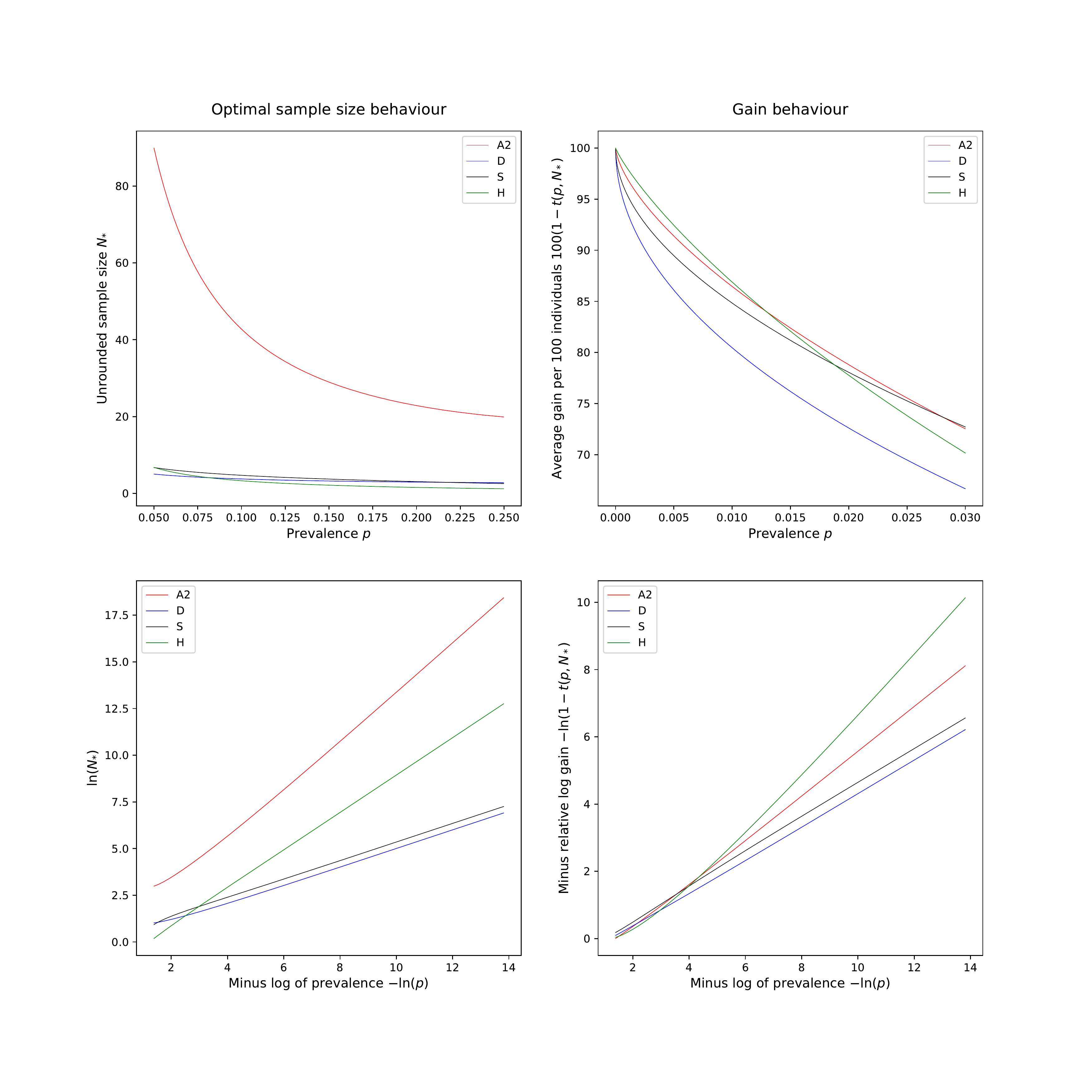}
\caption{Graph showing behaviour of optimal pool sizes and gains on original and log--log scales.} 
\label{fig:comparing_schemes1}
\end{figure}

\smallskip\noindent\emph{Halving scheme {H}}. $\rhd$ This scheme resembles divide and conquer sorting and can be described by the following steps.
\begin{itemize}
    \item[Step 1.] Test initial pooled cohort. If it tests negative, finish; if it tests positive, proceed to Step 2.
    \item[Step 2.] Divide the cohort into two approximately equal parts consisting of the first and second halves and apply the whole algorithm (starting from Step 1) to the two obtained parts recursively.
\end{itemize}
It is difficult to trace back the first reference discussing this scheme in detail. To our best knowledge, it was treated already by Johnson et. al. \cite{johnson_inspection_1991}. However, it appears that its asymptotic analysis was first accomplished not so long ago by Zamman and Pippenger \cite{zaman+pippenger:2016} whereas in \cite{Skorniakov_Leipus_2021} it was discussed a fresh without a strict focus on asymptotic regime when $p\to 0+$. There it was shown that rounded optimal configuration
\begin{equation*}
    N_{opt}^{(H)}\in\{\lfloor 1/(2\log_2(1/q))\rfloor,\lfloor 1/(2\log_2(1/q))\rfloor+1\}
\end{equation*}
leading thereby to the following asymptotic relationships:
\begin{gather*}
    N_*^{(H)}=-\frac{1}{2\log_2 q}(1+o(1)),\quad t^{(H)}(q,N_*^{(H)})=-({2p}\log_2 p)(1+o(1)),\\
    G^{(H)}(q)=(1+{2p}\log_2 p)(1+o(1)),\ p\to 0+.\ \lhd
\end{gather*}
Figure \ref{fig:comparing_schemes1} shows the behaviour of $N_*$ and gain $G$ for the case of A2 and the three schemes described above. Note that, in case of A2, total optimal pool size $N_*^{(A2)}=n_{min}^2$ with $n_{min}$ given in Theorem \ref{t:main}. Due to that raise to the square, $N_*^{(A2)}$ grows to infinity much faster than the counterparts of the remaining schemes, and, because of this, in the top left sub-figure, the range of $p$ starts quite far from the origin and the accompanying bottom left sub-figure on the log scale is given. It clearly depicts the relationships 
\begin{equation}\label{e:Ngrowth_on_log_scale}
    \ln\left(N_*^{(D)}\right), \ln\left(N_*^{(S)}\right) \sim -\frac{1}{2}\ln(p),
    \ln\left(N_*^{(H)}\right)\sim -\ln(p),\ln\left(N_*^{(A2)}\right)\sim -\frac{4}{3}\ln(p),
\end{equation}
following from the formulae given above and clearly showing that the asymptotic slope of $\ln\left(N_*^{(A2)}\right)$ on the $-\ln p$ scale is the largest one.

\begin{figure}[ht!]
\centering
\includegraphics[width=15cm, height=12cm]{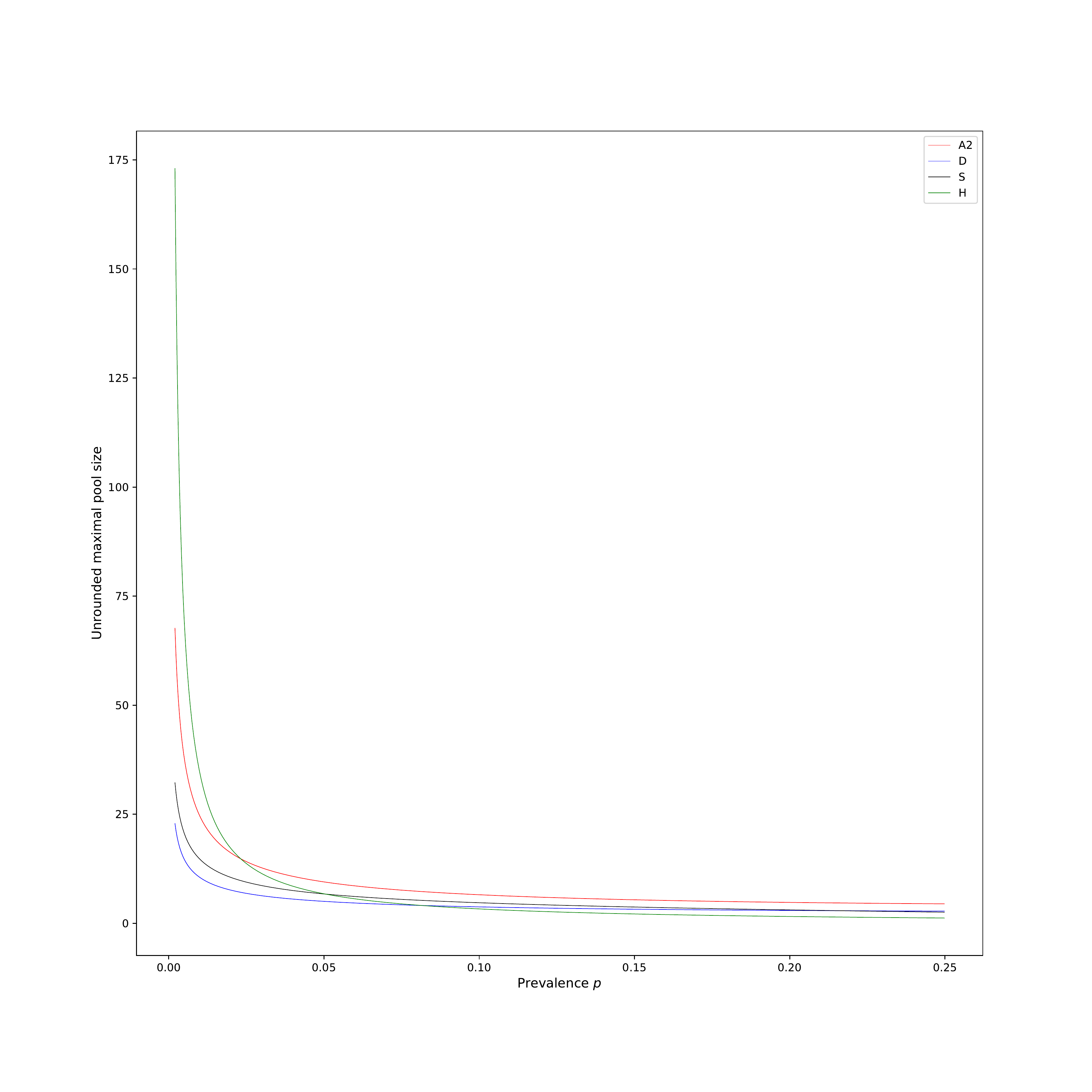}
\caption{Graph showing behaviour of maximal tested pool size for different schemes.} 
\label{fig:comp_max_pool}
\end{figure}

Turning to sub-figures on the right, one sees that there are ranges where A2 outperforms other schemes. Numerical solutions are as follows:
\begin{itemize}
    \item $R_{A2,D}=\{p\in(0,0.249790): G^{(A2)}(q)>G^{(D)}(q)\}=(0,0.115589)$;
    \item $R_{A2,S}=\{p\in(0,0.249790): G^{(A2)}(q)>G^{(S)}(q)\}=(0,0.028071)$;
    \item $R_{A2,H}=\{p\in(0,0.249790): G^{(A2)}(q)>G^{(H)}(q)\}=(0.012936,0.220788)$.
\end{itemize}
Within these ranges, numerically estimated maximal differences are
\begin{itemize}
    \item $\max_{p\in R_{A2,D}} 100(G^{(A2)}(q)-G^{(D)}(q))=6.2179$ at $p=0.017128$;
    \item $\max_{p\in R_{A2,S}}100(G^{(A2)}(q)-G^{(S)}(q))=1.9342$ at $p= 0.003984$;
    \item $\max_{p\in R_{A2,H}}100(G^{(A2)}(q)-G^{(H)}(q))=6.5951$ at $p=0.104908$.
\end{itemize}
Finally, comparing behaviour of optimal pool sizes one has that $\sqrt{N_*^{(A2)}}$ (and, therefore, $N_*^{(A2)}$ as well) exceeds $\max(N_*^{(D)},N_*^{(S)})$ for all $p\in(0,0.249790)$. For scheme H, however, the following holds true:
\begin{gather*}
    \sqrt{N_*^{(A2)}}<N_*^{(H)} \quad\text{for}\quad p\in(0,0.023178];\\
    \sqrt{N_*^{(A2)}}>N_*^{(H)} \quad\text{for}\quad p\in[0.023179,0.249790);\\
    {N_*^{(A2)}}>N_*^{(H)}\quad\text{for all}\quad p\in(0,0.249790).
\end{gather*}
Figure \ref{fig:comp_max_pool} provides a zoomed in visual illustration. Regarding the role and distinction between $\sqrt{N_*^{(A2)}}$ and ${N_*^{(A2)}}$, consult the discussion Section \ref{s:discussion}.

\subsubsection{Optimal configuration when the prevalence is unknown}\label{ss:example2}

In reference \cite{Malinovsky_note_on_minimax}, the authors sought for the pool size leading to optimal testing by making use of scheme D when the prevalence is unknown. In their work, two approaches were used. Both (approaches) were based on the following loss function. Given scheme $X$, define
\begin{equation}\label{e:base_loss}
    L^{(X)}(q,N)=t^{(X)}(q,N)-t^{(X)}(q,N_{opt}^{(X)}),
\end{equation}
where $N^{(X)}_{opt}=N^{(X)}_{opt}(q)=\argmin_{N\in\{1,2,\ldots\}}t^{(X)}(q,N)$ is the optimal configuration when the prevalence $p$ (and hence $q=1-p$) is known. It is clear that $L^{(X)}(q,N)\geq 0\ \forall (q,N)\in(0,1)\times\Nd$ and, for a given $q$, $L^{(X)}(q,N)=0$ precisely when $N=N_{opt}^{(X)}(q)$. 

In what follows, to distinguish between $N_{opt}^{(X)}(q)$ and optimal configuration suitable for unknown $q$'s, the latter configuration is denoted by $N_{\star}^{(X)}$.

The first approach in \cite{Malinovsky_note_on_minimax} was to make use of mini--max strategy and take $N_{\star}^{(X)}$ as a minimizer of
\begin{equation}\label{e:mini_max_func}
    \{1,2,\ldots\}\ni N\mapsto \sup_{q} L^{(X)}(q,N).
\end{equation}
The second approach was to make use of Bayesian paradigm and, after putting the prior $\pi$ on $q$, to take $N_{\star}^{(X)}$ as a minimizer of
\begin{multline}\label{e:Bayesian_func}
    \{1,2,\ldots\}\ni N\mapsto \Mean_{\pi} L^{(X)}(q,N)=\Mean_{\pi}t^{(X)}(q,N)-c(\pi),\quad\\ 
    c(\pi)=\Mean_{\pi}t^{(X)}(q,N_{opt}^{(X)}).
\end{multline}

In this example, we have adopted both approaches to the case of A2. When using the Bayesian one, $\pi$ was taken uniform over $(0.750210,1)$. Thereby, we have modeled situation when the only prior information is that application of A2 makes sense (see Corollary \ref{c:n_opt}). Also, we have modified \eqref{e:Bayesian_func} and used
\begin{equation}\label{e:Bayesian_func2}
    \{1,2,\ldots\}\ni N\mapsto \Mean_{\pi} (L^{(X)}(q,N))^2=\Mean_{\pi}\left(t^{(X)}(q,N)-t^{(X)}(q,N_{opt}^{(X)})\right)^2
\end{equation}
instead. To justify our choice, note that, in \eqref{e:Bayesian_func}, $c(\pi)$ does not depend on $N$. Therefore, minimization of the target function amounts to minimization of $\{1,2,\ldots\}\ni N\mapsto \Mean_\pi t^{(X)}(q,N)$. This way important information carrying function $q\mapsto t^{(X)}(q,N_{opt}^{(X)})$ remains unutilized. 

\begin{figure}[ht!]
\centering
\includegraphics[width=15cm, height=12cm]{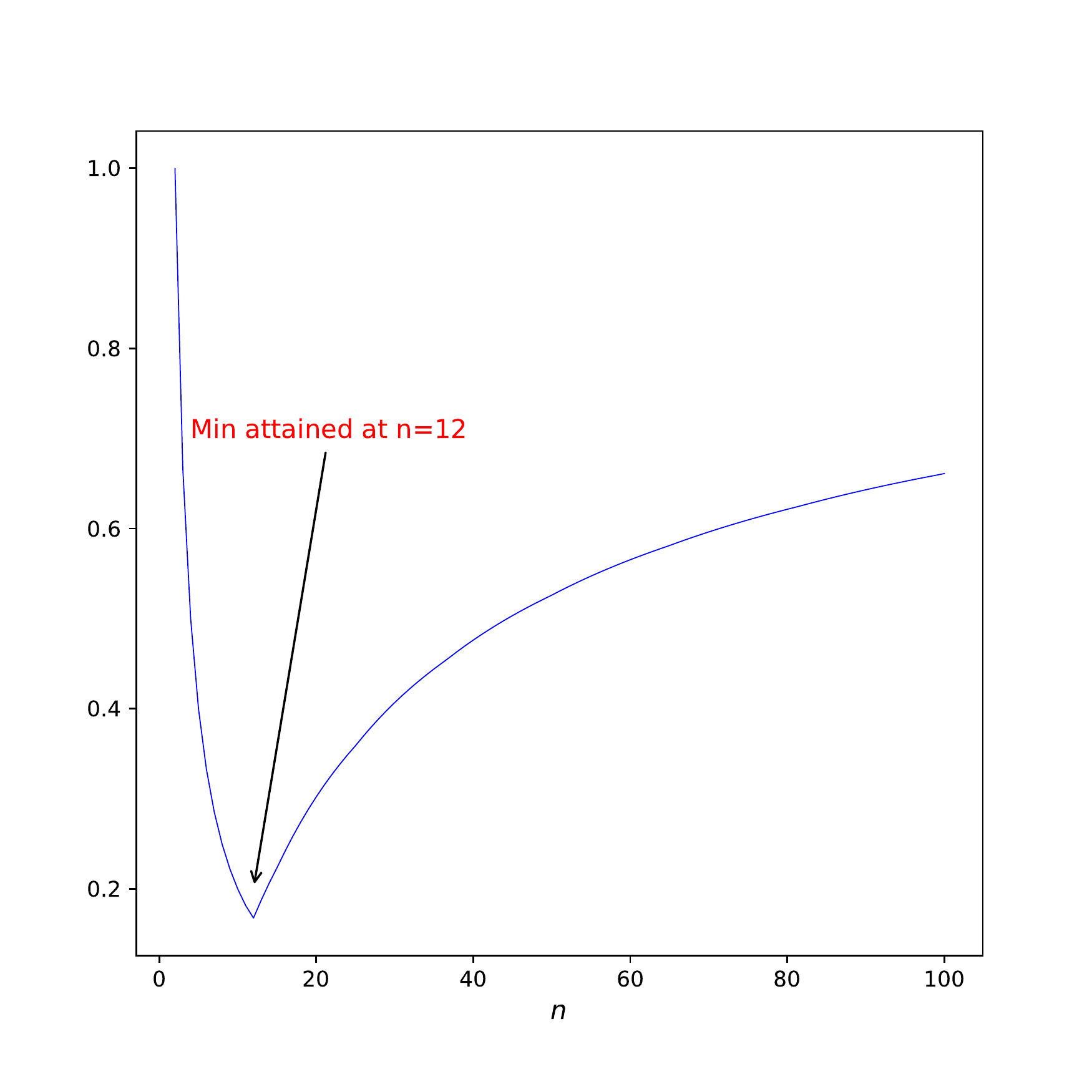}
\caption{Graph of $\{2,3,\ldots\}\ni n \mapsto \sup_{q} (t^{(A2)}(q,n)-t^{(A2)}(q,n_{opt}^{(A2)}))$  } \label{fig:minimax_graph}
\end{figure}

\begin{figure}[ht!]
\centering
\includegraphics[width=15cm, height=12cm]{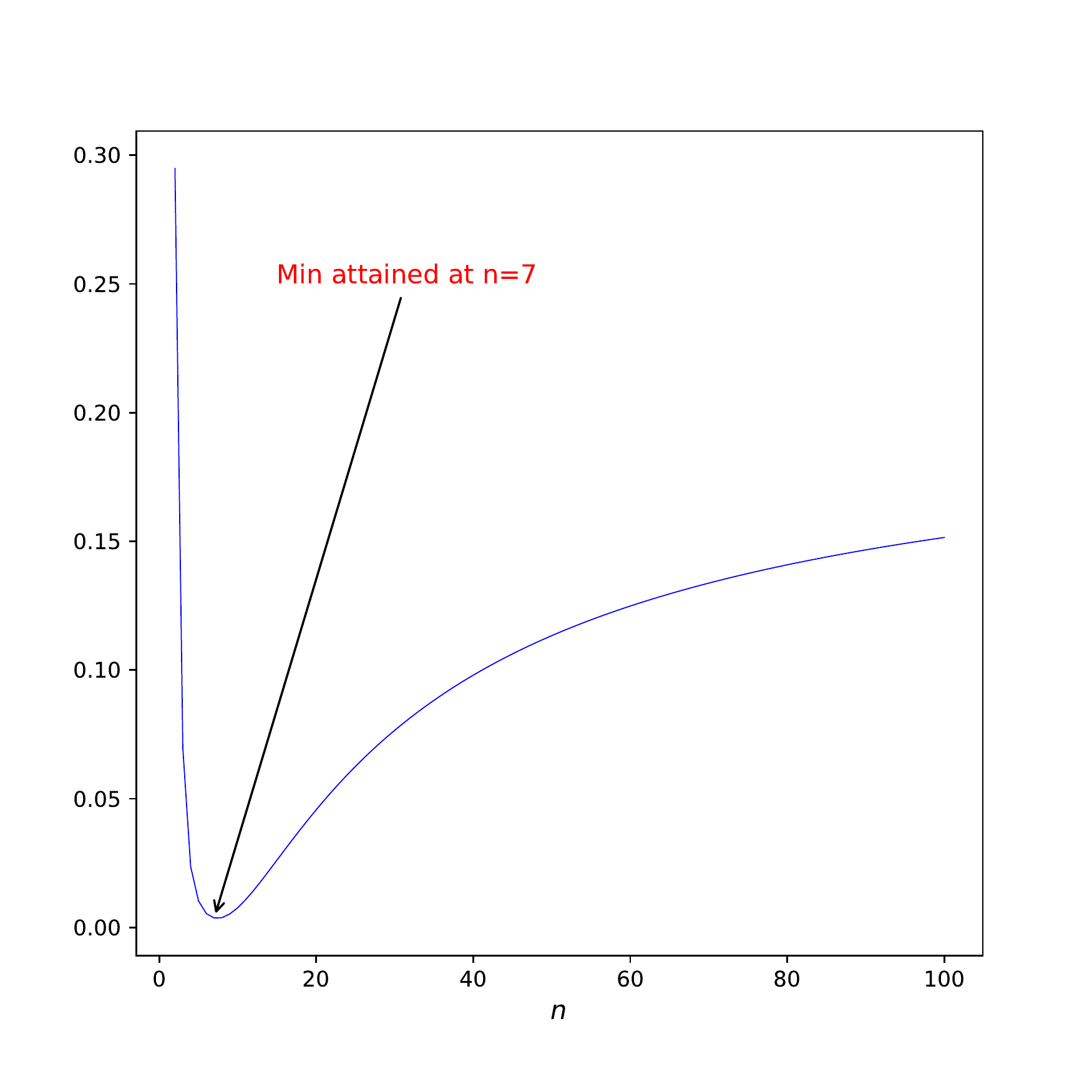}
\caption{Graph of $\{2,3,\ldots\}\ni n \mapsto \Mean_{\pi} (t^{(A2)}(q,n)-t^{(A2)}(q,n_{opt}^{(A2)}))^2$  } \label{fig:bayesian_graph}
\end{figure}

Figures \ref{fig:minimax_graph}--\ref{fig:bayesian_graph} show graphs of \eqref{e:mini_max_func} and \eqref{e:Bayesian_func2} for the case of $X=\mathrm{A2}$ and the previously mentioned prior $\pi$. Note that, adopting the above to our case, in \eqref{e:base_loss}, \eqref{e:mini_max_func}, \eqref{e:Bayesian_func}, we have used function $t$ defined by \eqref{e:t_q_n_A2}. Numerical estimation yielded the following values:
\begin{itemize}
    \item $N_{\star}^{(A2)}=12^2$ for the case of mini--max approach;
    \item $N_{\star}^{(A2)}=7^2$ for the case of Bayesian approach.
\end{itemize}
Finishing, it is important to note that, though the strategy discussed above leads to suboptimal testing in a stable environment where reliable estimation of prevalence is possible, it appears to be a reasonable strategy when the prevalence is varying rapidly and is difficult to capture by data at hand. Therefore, at least in the initial stage, it can be considered as a good alternative for optimal testing during pandemics like COVID--19.

\section{Discussion}\label{s:discussion}
It was already mentioned in the introductory Section \ref{s:intro} that, to our best knowledge, \cite{hudgens_optimal_2011} is the only reference where A2 was treated in the same way like we did here. Therefore, we first discuss our input in comparison with \cite{hudgens_optimal_2011} and then turn to a more general setting. 

\subsection{Comparison with the previous work}\label{ss:discussion1}

\begin{figure}[ht!]
\centering
\includegraphics[width=15cm, height=12cm]{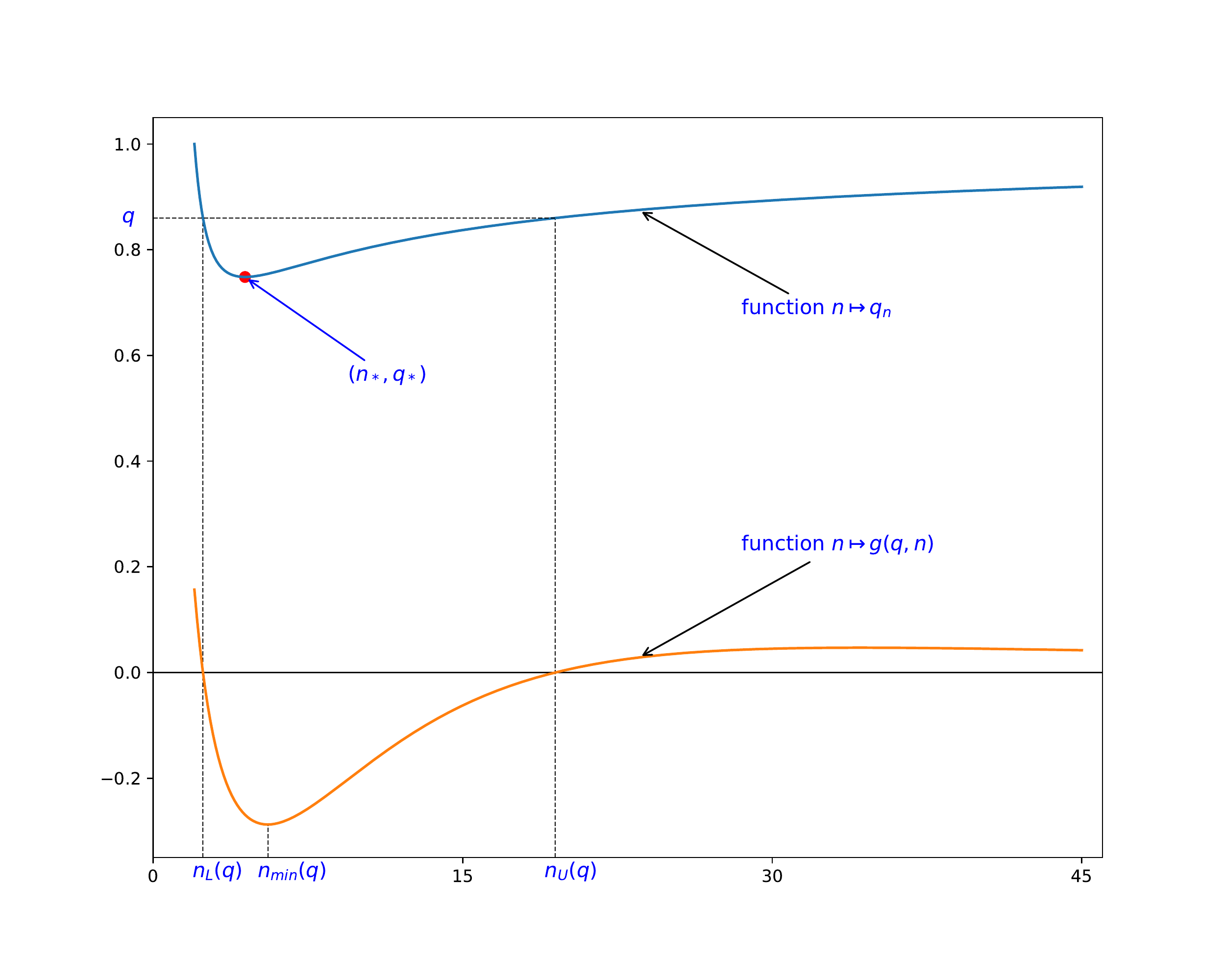}
\caption{Graph illustrating relationships of $n\mapsto q_n$ and related quantities. The  lower curve corresponds to $q=0.86$.} 
\label{fig:q_n_graph}
\end{figure}

Figure \ref{fig:q_n_graph} illustrates the behaviour of the function $n\mapsto q_n$ appearing in the proof of Theorem \ref{t:main}. The minimum of this function, denoted by $q_*$ is an exact lower bound of the region where A2 is efficient \emph{on a continuous scale} in a sense that, given prevalence $p\in(0,1-q_*)$ (or, alternatively, $q\in(q_*,1)$), the achievable minimum of $(2,\infty)\ni n\mapsto t^{(A2)}(q,n)$ is strictly smaller than $1$. This always holds true for the unrounded optimal configuration, i.e., the minimizer of $(2,\infty)\ni n\mapsto t^{(A2)}(q,n)$ on the continuous scale. In \cite{hudgens_optimal_2011}, the authors also provide a region of this kind. We utilize that region in Corollary \ref{c:n_opt} and denote it $(q_5,1)$. Remark \ref{r:about_region_truncation} explains why it can be thought of as a region where A2 makes sense from the practical point of view and how this can be derived from our results. Linking our work to theirs, it is important to mention that, deriving their proofs, Hudgens and Kim \cite{hudgens_optimal_2011} operated on the discrete scale and obtained this region as the one where $\{5,6,\ldots\}\ni n\mapsto t^{(A2)}(q,n)$ is efficient. Having proved that $n=2,3,4$ are never optimal, they have also verified that $(q_5,1)$ is the region where practical application of A2 makes sense.

Summing up, in this direction, our input adds the missing part to the theoretical characterization of A2 yet does not bring novelty from the applied point of view. Turning to formulae \eqref{e:A2_n_optimal}, \eqref{e:set_of_optimal_vals}, the situation is essentially different. Kim and Hudgens \cite{hudgens_optimal_2011} have given lower and upper bounds on the configuration $\sqrt{N_{opt}^{(A2)}}=n_{opt}^{(A2)}\in\{5,6,\ldots\}$ and ascertained that the optimal configuration is unique. These bounds facilitate numerical computations and allow some insights on the analytical scale yet they do not accommodate all the benefits brought by formulae \eqref{e:A2_n_optimal}--\eqref{e:set_of_optimal_vals}:
\begin{enumerate}
    \item \eqref{e:set_of_optimal_vals} provides \emph{a ready to use} expressions for $N_{opt}^{(A2)}$ which are practically important when it comes to $p$'s ranging near the origin where numerical solutions might cause computational problems;
    \item \eqref{e:A2_n_optimal} and \eqref{e:set_of_optimal_vals} enable to obtain analytical comparisons of A2 with other schemes. An example of such analysis is given in Subsection \ref{ss:example1} when deriving relationships \eqref{e:Ngrowth_on_log_scale}. Without knowing analytical expressions of $N_{opt}^{(A2)}$ it would be impossible to contrast the asymptotic slope of $-\ln p\mapsto\ln N_*^{(A2)}(p)$ to the analogous asymptotic slopes of other schemes considered there. One can go even further. Namely, if $X$ is another scheme of interest, one can compare analytically the behaviour of $t^{(A2)}(q,N_{opt}^{(A2)})$ with $t^{(X)}(p,N_{opt}^{(X)})$ as $p\to 0+$ (see a comment in the forthcoming Subsection \ref{ss:discussion2} regarding an importance of analysis of this kind). It is again worth to stress up that, in this asymptotic setting, validity of pure numerical analysis is always questionable whereas analytical result leads to accurate and definite analysis.
\end{enumerate}

\subsection{Discussion of other aspects}\label{ss:discussion2}
Several authors have already noted that array based algorithms can be more efficient in certain settings
\cite{Berger-2000}, \cite{kim_comparison_2007}, \cite{kim_three-dimensional_2009}. Our findings confirm these observations: comparisons given in Subsection \ref{ss:example1} demonstrated that there were regions where A2 performed better than other considered schemes. The same behaviour is expected for other unconsidered schemes satisfying binomial testing assumptions. Therefore, though initially applied in a frame of genetic a screening \cite{barillot_theoretical_1991}, \cite{amemiya_two-dimensional_1992}, \cite{bruno_efficient_1995}, A2 is an appropriate candidate for other GT applications as well. Since GT is most effective when the prevalence $p$ is low, for any scheme $X$, asymptotic behaviour of $t^{(X)}(p,N_{opt}^{(X)})$ as $p\to 0+$, is an important characteristic. Among the schemes contrasted in Subsection \ref{ss:example1}, A2 takes an intermediate position (with respect to this characteristic) because\footnote{the relationships are not difficult to justify by making use of the formulae \eqref{e:set_of_optimal_vals}}
\begin{gather*}
    \lim\limits_{p\to 0+}\frac{t^{(A2)}(q,N_{opt}^{(A2)})}{t^{(X)}(q,N_{opt}^{(X)})}=0\quad\text{for}\quad X\in\{D,S\},\\ \nonumber
    \text{and}\quad
    \lim\limits_{p\to 0+}\frac{t^{(H)}(q,N_{opt}^{(H)})}{t^{(A2)}(q,N_{opt}^{(A2)})}=0.\label{e:t_rates}
\end{gather*}
This, however, comes at a cost of a quickly increasing \emph{maximal tested pool size} at optimal configuration with relationships inverse to those above:
\begin{gather}
    \lim\limits_{p\to 0+}\frac{N_{opt}^{(X)}(q)}{\sqrt{N^{(A2)}_{opt}(q)}}=0\quad\text{for}\quad X\in\{D,S\},\label{e:N_rates1}\\
    \text{and}\quad
    \lim\limits_{p\to 0+}\frac{\sqrt{N_{opt}^{(A2)}(q)}}{N_{opt}^{(H)}(q)}=0.\label{e:N_rates2}
\end{gather}
We would like to emphasize that, in case of A2, $\sqrt{N_{opt}^{(A2)}}=n_{opt}^{(A2)}$ (not $N_{opt}^{(A2)}$) is the maximal tested pool size (corresponding to the test of a distinct row/column) whereas in case of all the rest schemes it is $N_{opt}^{(X)},X\in\{D,H,S\}$. Thus the appearance of $\sqrt{N_{opt}^{(A2)}}$ in formulae \eqref{e:N_rates1}--\eqref{e:N_rates2} above. $N_{opt}^{(A2)}$ is the batch size needed to organize optimal testing since  $N_{opt}^{(A2)}=n_{opt}^{(A2)}$ samples never get tested as a single specimen. This distinction is important since in the applied setting one has to take into account the fact that the operating characteristics of the test might depend on the size of the tested pool. That is, sensitivity and/or specificity of the test might become unacceptably low when the pooled sample is formed out of $N\geq N_{critical}$ individual samples. Also, subjects might be dependent and it might be reasonable to take into account imperfectness of the test even on the single individual basis. These facts limit usage of many schemes satisfying binomial testing assumptions. E.g., scheme H, also known under the name of optimal testing algorithm and giving best asymptotic performance as $p\to0+$, requires large pool sizes. Inspecting figure \ref{fig:comp_max_pool}, one sees that A2 also exhibits such behaviour. Nonetheless, investigations giving theoretical characterization of the GT scheme considered are important because of the following reasons.

\begin{enumerate}
    \item[Convenience.] When choosing between two realizable schemes, the preference does not always fall on the optimal one. If it turns out that the optimal scheme yields a small surplus in gain, it can be exchanged to a more simply realizable competitor. E.g., scheme D was recently employed in Lithuania \cite{PoolLT} for massive COVID-19 testing and a large number of other countries do the same \cite{Wiki_2020}. Turning to A2, it appears that one of the reasons of its emergence in genetic applications was an operational convenience.
    \item[Tolerable errors.] Though there exist a lot of generalizations allowing imperfectness of the test (e.g., \cite{kim_comparison_2007}, \cite{kim_three-dimensional_2009}, \cite{habtesllassie_array-based_2015}), as noted by several authors \cite{Berger-2000}, \cite{hudgens_optimal_2011}, \cite{malinovsky_revisiting_2019}, the schemes assuming perfect tests can be quite accurate since modern tests exhibit very small errors.
    \item[Benchmarking.] BTA based schemes, being more simple to treat analitically, provide theoretically justified benchmark thresholds for more elaborated schemes which assume imperfectness of the test and/or other specific conditions (e.g., testing outcome dependence on subject specific characteristics).
    \item[Basement.] More elaborated schemes may emerge on the basis of the simpler ones. An example is the scheme S introduced in Subsection \ref{ss:example1}. It was built on the top of scheme D. In \cite{kim_comparison_2007} and \cite{habtesllassie_array-based_2015}, A2 serves as a basis for extensions incorporating imperfectness of the test and dilution effect.
\end{enumerate}

\bibliographystyle{plain}

\appendix
\section{Proofs}\label{a:proofs}
\begin{rem}
In the proof of Theorem \ref{t:main}, Step 1 is a repetition of Lemma 1 in \cite{hudgens_optimal_2011}. We have decided to rewrite it here because it is very short and, along the way, some notions used in the sequel appear.
\end{rem}

\emph{Proof of Theorem \ref{t:main}.}
\emph{Step 1.} We first show that, for any fixed $n\in (2,\infty)$, there exists a unique $q_n\in(0,1)$  such that
\begin{equation*}
    g(q_n,n)=0, \quad g(q,n)<0 \ \forall q\in(q_n,1)\quad\text{and}\quad g(q,n)>0\ \forall q\in (0,q_n).
\end{equation*}
To this end, note that 
\begin{multline*}
    \frac{\partial}{\partial q} g(q,n)=-2nq^{n-1} + (2n-1)q^{2n-2}=\\
    -2nq^{n-1}\left(1-\frac{2n-1}{2n}q^{n-1}\right)<0\ \forall n\in(2,\infty).
\end{multline*}
Thus, given $n\in(2,\infty)$, $q\mapsto g(q,n)$ is decreasing on $(0,1)$. Since $g(0+,n)=\frac{2}{n}>0$ and $g(1-,n)=\frac{2}{n}-1<0$, the claim holds true. 

\emph{Step 2.} From \emph{Step 1} it follows that, for any fixed $q\in (0,1)$, we have a well defined function $n\mapsto q_n$ which is given implicitly by equation $g(q_n,n)=0$. Since this function is continuous, its range $I\subset(0,1)$ is an interval. Further, for $\eps\in(0,1)$ and $n=2+\eps$,
\begin{multline*}
    g(1-\eps,n)=\frac{2}{n} - 2(1-n\eps+O(\eps^2))+(1-(2n-1)\eps+O(\eps^2))=\\
    \frac{-\eps}{2+\eps}+\eps + O(\eps^2)>0, 
\end{multline*}
provided $\eps$ is small enough. Hence, analysis accomplished in \emph{Step 1} implies that $q_{2+\eps}\in (1-\eps,1)$. Therefore, $I =(q_*,1)$ for some $q_*\in(1/2,1)$. To justify the lower bound $1/2$, note that
\begin{equation*}
    g(1/2,n)>0 \Longleftrightarrow 2^n>n\left(1-\frac{1}{2^n}\right).
\end{equation*}
Since the right hand side holds true for all $n>2$, it follows that $g(1/2,n)>0\ \forall n> 2$. Also note that $g(q,n)>0\ \forall\,(q,n)\in(0,q_*)\times (2,\infty)$ since an opposite contradicts the definition of $q_*$.

\emph{Step 3.} Fix $q\in(q_*,1)$. By \emph{Step 1--Step 2}, $g(q,n)=0$ has at least one solution $n=n(q)$ (suffices it to take $n$ such that $q_n=q$). To show that there are two solutions, put $c=\frac{1}{2q}$, make a change of variable $q^n=x$, and rewrite $g(q,n)=0$ in a form
\begin{equation}\label{e:g_eq_0_first}
    -\ln q=-\ln x(x-cx^2).
\end{equation}
Consider function $h(x)=-(x-cx^2)\ln x$ for $x\in(0,1)$. Note that
\begin{equation}\label{e:hderiv}
    \frac{\d}{\d x}h(x)=-\left((1-cx)+(1-2cx)\ln x\right)=0\Longleftrightarrow
    -\ln x=\frac{1-cx}{1-2cx}.
\end{equation}
Since $1-cx>1-c>0$ and $-\ln x>0$ for all $x\in(0,1)$, it follows that $1-2cx>0$ as well, provided $x$ solves \eqref{e:hderiv}. Therefore, the range of possible solutions of \eqref{e:hderiv} shrinks to $(0,q)$. Moreover, relationships $\lim\limits_{x\to 0+}\frac{\d}{\d x}h(x)=\infty$, $\lim\limits_{x\to 1-}\frac{\d }{\d x}h(x)=c-1<0$ imply that \eqref{e:hderiv} has at least one solution. Since $\frac{1-cx}{1-2cx}+\ln x= 1+ \frac{cx}{1-2cx}+\ln x$ increases on $(0,q)$, the solution is unique. Denote it $x_0$. Based on the sign of the derivative, we have that $h\uparrow$ on $(0,x_0)$ and $h\downarrow$ on $(x_0,1)$. Hence, at $x_0$, $h$ attains its maximum and \eqref{e:g_eq_0_first} admits exactly two solutions if $h(x_0)>-\ln q$, one solution if $h(x_0)=-\ln q$, and has no solutions if $h(x_0)<-\ln q$. By the choice of $q$ (recall that $q>q_*$), the last case can not hold. To exclude the second one, note that the function $\left(\frac{1}{2},\frac{1}{2q_*}\right)\ni c\mapsto x_0(c)$ is well defined and decreasing since
\begin{multline*}
    -\ln x_0=\frac{1-cx_0}{1-2cx_0}\imply\\
    -\frac{\d}{\d c}\ln x_0(c)=-\frac{\frac{\d}{\d c}x_0(c)}{x_0(c)}=\frac{\d}{\d c}\frac{1-cx_0}{1-2cx_0}=\frac{x_0+c\frac{\d}{\d c}x_0(c)}{(1-2cx_0)^2}\imply\\ \frac{\d}{\d c}x_0(c)=\frac{-x_0^2}{(1-2cx_0)^2+cx_0}<0.
\end{multline*}
Therefore, $(q_*,1)\ni q\mapsto x_0(q)$ is increasing. Taking into account that $q\mapsto -\ln q$ is decreasing, we finally deduce that $h(x_0)>-\ln q$ for all $q\in(q_*,1)$. The monotonicity of $q\mapsto x_0(q)$ and $q\mapsto -\ln q$ also leads to conclusion that $q_*$ can be solved from equation
\begin{equation*}
    -\ln q_*=h(x_0(q_*))
\end{equation*}
along with a unique $n_*\in(2,\infty)$. Hence {(i)}.
    
\emph{Step 4.} Assume the setting of \emph{Step 3}. Let $0<x_L=q^{n_U}<x_U=q^{n_L}<q$ denote two solutions of \eqref{e:g_eq_0_first}. By above, $h(x)>-\ln q\Longleftrightarrow x\in (x_L,x_U)$. Reverting to $(0,\infty)\ni n\mapsto g(q,n)$, this reads as $g(q,n)<0\Longleftrightarrow n\in (n_L,n_U)$. Note that $n\mapsto g(q,n)>0$ in the neighborhood of $\infty$. Also, from \emph{Step 1}--\emph{Step 2}, we have that $n\mapsto g(q,n)>0$ in the right neighborhood of 2 and that $n_U\in (2,\infty)$. Therefore, continuity of $n\mapsto g(q,n)>0$ yields that $n_L\in(2,\infty)$ as well. Finally, it is clear that $n_L<n_*<n_U$ (see figure \ref{fig:q_n_graph} for a graphical illustration). Hence (ii).

\emph{Step 5.} In this step, we identify the number and location of zeroes of the derivative of $(2,\infty)\ni n\mapsto g(q,n)$ having fixed $q\in (q_*,1)$. From analysis given in \emph{Step 4}, it follows that $(2,\infty)\ni n\mapsto g(q,n)$ has at least two extremes: there must be a minimum in $(n_L,n_U)$ (since function is negative here), and maximum in $(n_U,\infty)$ (since the function is positive here and $\lim_{n\to\infty} g(q,n)=0+$. To see that there are no other extremes, consider equation
\begin{equation}\label{e:g_deriv_raw}
    \frac{\partial}{\partial n}g(q,n)=-\frac{2}{n^2}-2q^n\ln q+2q^{2n-1}\ln q=0,   
\end{equation}
and rewrite it by making use of notions introduced in \emph{Step 3} as follows:
\begin{equation}\label{e:g_deriv}
    -x\ln x=\sqrt{-\ln q}\sqrt{\frac{x}{1-2cx}}, x\in \left(0,\frac{1}{2c}\right)=(0,q).    
\end{equation}
Next, note that: 
\begin{itemize}
    \item $h_1(x)=-x\ln x$ is strictly convex--up and positive on $(0,1)$ with $\lim\limits_{x\to 0+}h_1(x)=\lim\limits_{x\to 1-}h_1(x)=0+$;
    \item $h_2(x)=\sqrt{\frac{x}{1-2cx}}$ is strictly positive and increasing on $\left(0,\frac{1}{2c}\right)=(0,q)$, it has one inflection point and $\lim\limits_{x\to 0+}h_2(x)=0$, $\lim\limits_{x\to q-}h_2(x)=\infty$.
\end{itemize}
Taking this information into account, we conclude that \eqref{e:g_deriv} can have at most two solutions and confirm thereby the assertion stated above.

\emph{Step 6.} It remains to justify expression \eqref{e:A2_n_optimal}. Let 
\begin{equation}\label{e:nt_def}
    n(q,t)= \frac{1}{p^{\frac{2}{3}}} + \frac{1}{2p^{\frac{1}{3}}} + 0.2 + 3p^2 + t,\quad t\in [0,1].
\end{equation}
It suffices to prove that, for all $q\in[0.755,1)$, the following statements hold true:
\begin{itemize}
    \item[(a)] $\max(g(q,n(q,0)),g(q,n(q,1))<0$; and
    \item[(b)] $\exists t\in[0,1]:\frac{\partial}{\partial n}g(q,n)\Big|_{n=n(q,t)}=0$.
\end{itemize}
Analytical calculations behind (a) and (b) are standard yet very lengthy and tedious. Therefore, we omit the details and end up with a graphical proof and a sketch of the analytical one.

\begin{figure}[ht!]
\centering
\includegraphics[width=15cm, height=12cm]{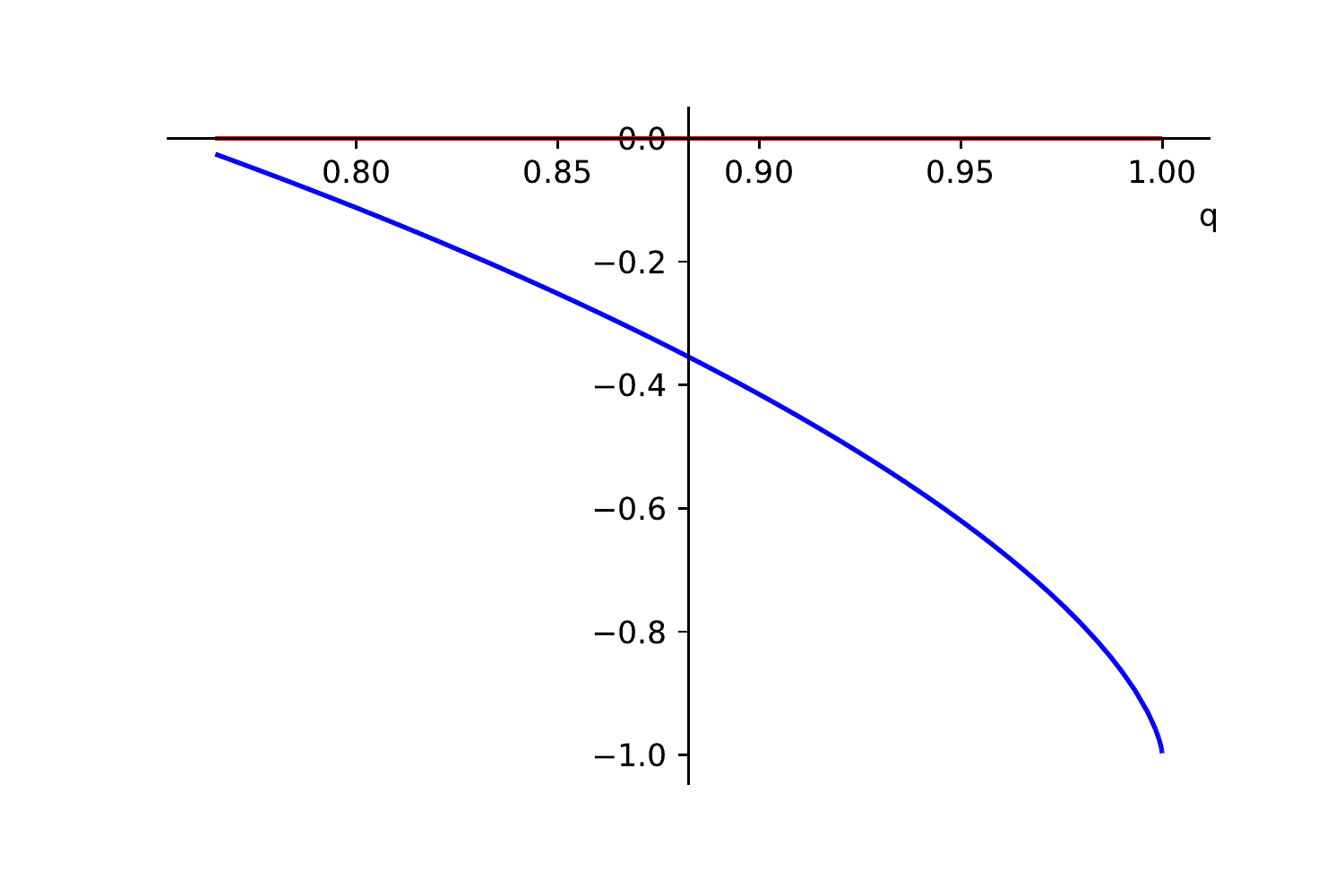}
\caption{Graph of $q\mapsto\max(g(q,n(q,0)),g(q,n(q,1)))$ for $ {q\in[0.765,1)}$. For reference, function identically equal to 0 is plotted in red.} \label{fig:g_01_max1}
\end{figure}

\begin{figure}[ht!]
\centering
\includegraphics[width=15cm, height=12cm]{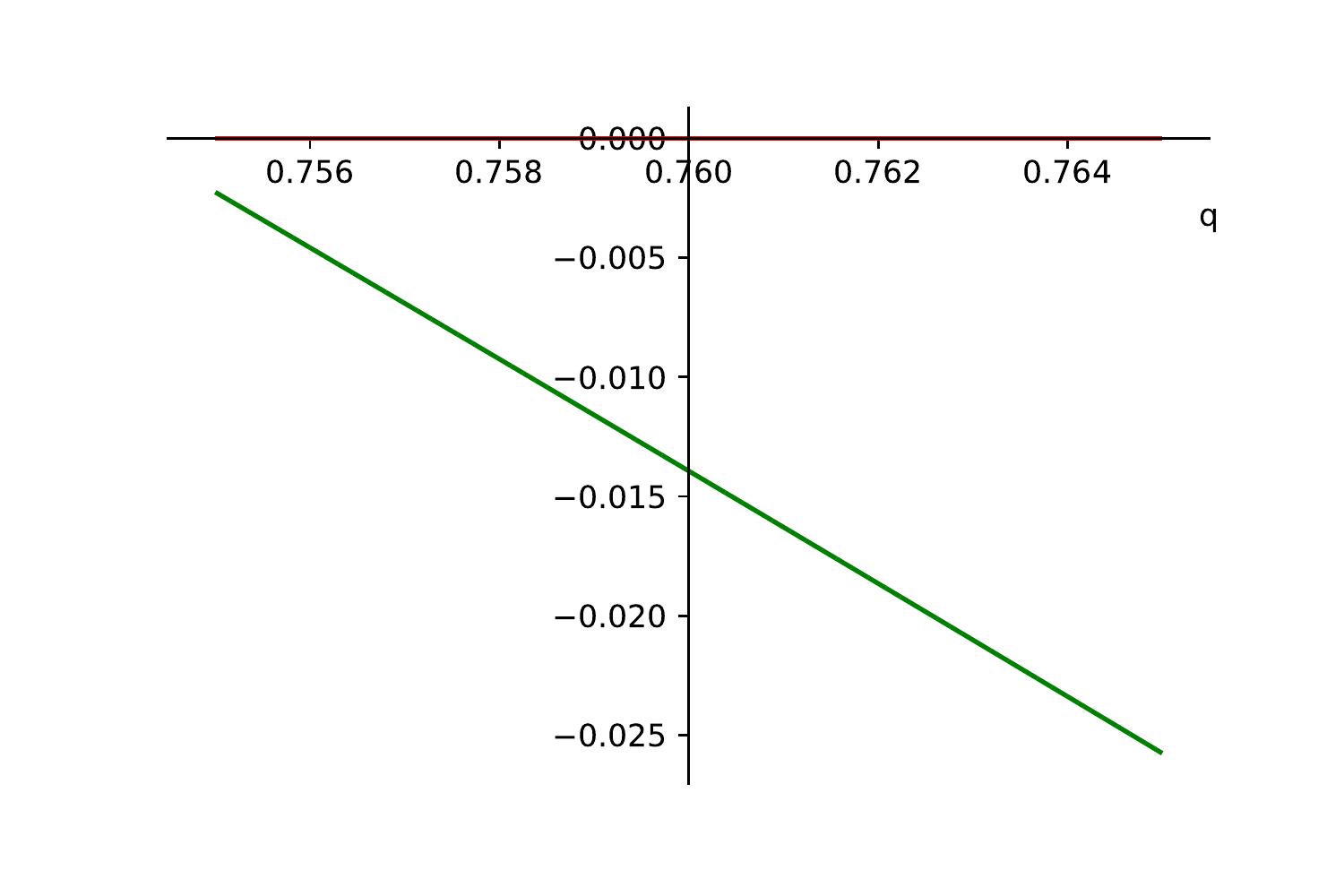}
\caption{Graph of $q\mapsto\max(g(q,n(q,0)),g(q,n(q,1))$ for ${q\in[0.755,0.765]}$. For reference, function identically equal to 0 is plotted in red.} \label{fig:g_01_max2}
\end{figure}

\begin{figure}[ht!]
\centering
\includegraphics[width=15cm, height=12cm]{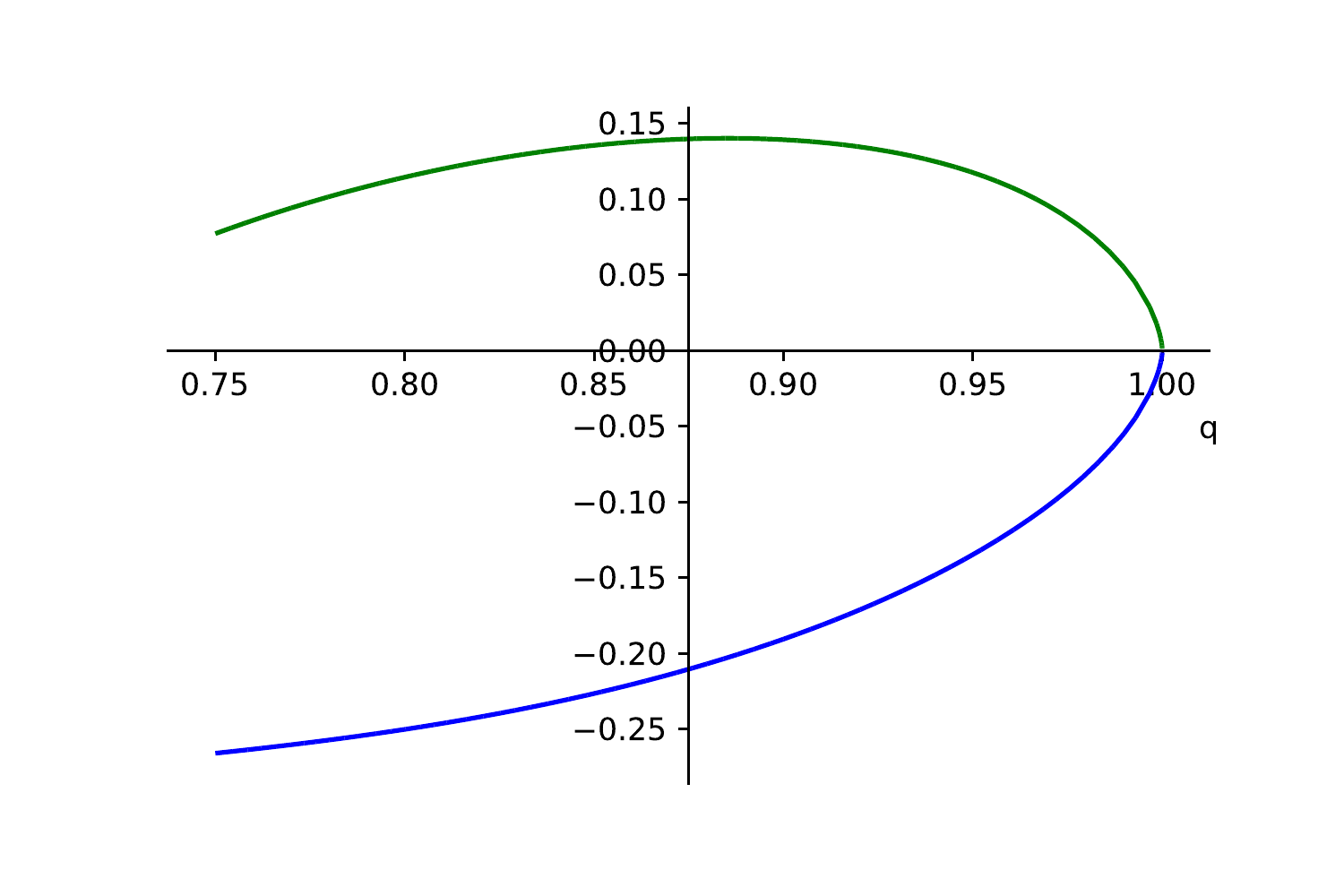}
\caption{Graphs of $q\mapsto h(0,q), q\in[0.755,1)$ (plotted in blue) and $q\mapsto h(1,q), q\in[0.755,1)$ (plotted in green).} \label{fig:h_01}
\end{figure}

Figures \ref{fig:g_01_max1}--\ref{fig:g_01_max2} show graph of $q\mapsto\max(g(q,n(q,0)),g(q,n(q,1)))$ from which it is evident that (a) holds. Analytical proof consists of the following steps.
\begin{itemize}
    \item[(s1)] Calculate $\frac{\partial}{\partial q}g(q,n(q,i)),i=0,1$.
    \item[(s2)] Check that $\frac{\partial}{\partial q}g(q,n(q,i))<0, i=0,1$ on $[0.755,1)$ and deduce that $g(q,n(q,i))$ decrease on [0.755,1). 
    \item[(s3)] Conclude that (a) indeed holds since 
    \begin{equation*}
    g(0.755,n(0.755,0))\approx-0.002258\quad\text{and}\quad g(0.755,n(0.755,1))\approx-0.013690.    
    \end{equation*}
    
\end{itemize}
Turning to (b), first rewrite \eqref{e:g_deriv_raw} as follows:
\begin{equation*}
    -{n^2}q^n\ln q(1-q^{n-1})=1.
\end{equation*}
Next, consider function $h(t,q)=-{n^2(q,t)}q^{n(q,t)}\ln q(1-q^{n(q,t)-1})-1$ with $n(q,t)$ given by \eqref{e:nt_def} and $q\in[0.755,1)$. Since $t\mapsto h(t,q)$ is continuous, it suffices to show that, for any $q\in[0.755,1)$, $h(0,q)<0$ and $h(1,q)>0$. Figure \ref{fig:h_01} shows graphs of $q\mapsto h(0,q),q\mapsto h(1,q)$. These confirm (b). Considering analytical part, the following is the suggested route.
\begin{itemize}
    \item[(s1)] Calculate $\frac{\partial^2}{\partial q^2}h(i,q),i=0,1$.
    \item[(s2)] Check that $\frac{\partial^2}{\partial q^2}h(0,q)>0$ whereas $\frac{\partial^2}{\partial q^2}h(1,q)<0$ on $[0.755,1)$ and deduce that $h(0,q)$ is convex downwards whereas $ h(1,q)$ is convex upwards on [0.755,1).
    \item[(s3)] By making use of Taylor's expansion, check that $\lim\limits_{q\to 1-}h(0,q)=0-$ and
    $\lim\limits_{q\to 1-}h(1,q)=0+$.
    \item[(s4)] Conclude that (b) indeed holds since 
    \begin{equation*}
        h(0,0.755)\approx-0.2645889 \quad \text{and}\quad h(1,0.755)\approx 0.081749.
    \end{equation*}
\end{itemize} \qed

\emph{Proof of Corollary \ref{c:n_opt}.} Uniqueness of $q_5$ was established in Step 1 of the proof of Theorem \ref{t:main}. Hudgens and Kim \cite{hudgens_optimal_2011} (Lemmas 2, 7, and 14) have demonstrated that $n_{opt}(q)\not \in\{2,3,4\}\ \forall q\in(0,1)$. From their results we also have that $n_{opt}(q)=5\ \forall q\in(q_*,0.755]$. It is straightforward to verify that 
\begin{equation*}
    \forall\ q\in (q_5,0.755]\ \left\lceil \frac{1}{p^{\frac{2}{3}}} + \frac{1}{2p^{\frac{1}{3}}} + 3p^2 + 1.2\right\rceil=5.    
\end{equation*}
Hence the claim in the region $(q_5,0.755]$. For $q\in[0.755,1)$, it follows from Theorem \ref{t:main} by noting that at least one of numbers in the set \eqref{e:set_of_optimal_vals} belongs to $\{n\in(2,\infty):t(q,n)<1\}$ because of (a) in Step 6 of proof of Theorem \ref{t:main}. \qed

\section{Tables and Figures}\label{app:tbl_and_figures}

The material contained in this appendix was produced by making use of an open source computer algebra system
SymPy \cite{sympy} and the python packages constituting the core kit for scientific numerical programming with python: NumPy \cite{numpy}, SciPy \cite{scipy}, Matplotlib \cite{matplotlib} and Pandas \cite{pandas_mckinney}, \cite{pandas}. The table accompanying Subsection \ref{ss:example1} is placed in our
\href{https://osf.io/bdxca}{Open Science Framework repository} (file A2\_table.csv). For the prevalence $p$ ranging over $(0,0.249790)$ at a step equal to $10^{-6}$ and for each scheme considered there, the following information is given: optimal (rounded) pool size $N_{opt}(q)$ and gain $G(q)=1-t(q,N_{opt}(q))$. In case of A2, $n^{(A2)}_{opt}(q)$ and $t^{(A2)}(q,n^{(A2)}_{opt}(q))$ are reported (see Subsection \ref{ss:example1}, Remark \ref{r:reparametrization} for the details). 

\end{document}